\documentclass[12pt]{amsart}
\usepackage{amssymb}
\newtheorem{theorem}{Theorem}

\newtheorem{corollary}[theorem]{Corollary}

\newtheorem{example}[theorem]{Example}

\newtheorem{proposition}[theorem]{Proposition}

\def \bR{\mathbb R}
\def \bC{\mathbb C}

\def \bS{\mathbb S}
\def \bN{\mathbb N}

\def \Oka{\mathcal O}

\begin{document}

\title[On the mixed Cauchy problem with data on singular conics]{On the
mixed Cauchy problem
 with data on singular conics}
\author[Peter Ebenfelt and Hermann Render]{Peter Ebenfelt and Hermann
Render}


\thanks{The first author is
supported in part by DMS-0401215.
The second author is
supported in part by Grant BFM2003-06335-C03-03 of the D.G.I. of Spain.}

\address{P. Ebenfelt: Department of Mathematics, University of California,
San Diego, La Jolla, CA 92093--0112, USA.}
\email{pebenfelt@math.ucsd.edu}
\address{H. Render: Departamento de Matem\'{a}ticas y Computaci\'{o}n,
Universidad de La Rioja, Edificio Vives, Luis de Ulloa s/n., 26004
Logro\~{n}o, Espa\~{n}a.} \email{render@gmx.de}

\maketitle

\begin{abstract} We consider a problem of mixed Cauchy type for
certain holomorphic partial differential operators whose principal
part $Q_{2p}(D)$ essentially is the (complex) Laplace operator to a
power, $\Delta^p$. We pose inital data on a singular conic divisor
given by $P=0$, where $P$ is a homogeneous polynomial of degree
$2p$. We show that this problem is uniquely solvable if the
polynomial $P$ is elliptic, in a certain sense, with respect to the
principal part $Q_{2p}(D)$.
\end{abstract}

\section{Introduction}

In this paper, we shall consider the mixed Cauchy problem for
holomorphic partial differential operators of the type
\begin{equation}
Lu= Q_{k}\left( D\right) u+\sum_{\left| \alpha \right| \leq
k_{0}}a_{\alpha }\left( z\right) D^{\alpha }u  \label{defLL}
\end{equation}
where $Q_k(D)$ is a non-trivial homogeneous, constant coefficient
partial differential operator of order $k$, the $a_{\alpha }(z)$ are
holomorphic functions in a domain $\Omega\subset\mathbb{C}^{n}$
containing $0$, and $k_{0}$ is a natural number $<k.$ We use
standard multi-index notation with $\alpha =\left( \alpha
_{1},...,\alpha _{n}\right) \in \mathbb{N}_{0}^{n}$,
$|\alpha|=\alpha_1+\ldots+\alpha_n$, and $D^{\alpha }$ denoting the
differential operator $$ D^\alpha:=\frac{\partial ^{\alpha
}}{\partial z^{\alpha }}=\frac{\partial ^{\left| \alpha \right|
}}{\partial z_{1}^{\alpha _{1}}...\partial z_{n}^{\alpha _{n}} }.$$
The principal symbol of the partial differential operator $L$ in
\eqref{defLL} is the homogeneous polynomial $Q_{k}\left(
\zeta\right) $ of degree $k$.

Let $\Oka(U)$ denote the space of holomorphic functions in $U\subset
\bC^n$. Clearly, $L$ defines a continuous linear operator $L\colon
\Oka(U)\to \Oka(U)$ for every $U\subset \Omega$. In general, this
linear operator is not injective. Indeed, if $L$ has, e.g., constant
coefficients, then it is well known (and easy to prove using the
idea of Fischer duality as in \cite{Shap89}; see also \cite{EbSh96})
that there are non-trivial entire solutions of $Lu=0$ as long as
$k\geq 1$ and, hence, $L\colon \Oka(U)\to \Oka(U)$ is never
injective in this case. The surjectivity of $L\colon \Oka(U)\to
\Oka(U)$ is more subtle, even when $L$ has constant coefficients,
and depends on the geometry of the domain $U$. The reader is
referred to \cite{HormC} for further information about this question
(see also \cite{APS}). However, it is an immediate consequence of
the classical Cauchy-Kowalevsky theorem that, for any domain $0\in
U\subset \Omega$, there exists a subdomain $0\in U'\subset U$ such
that the equation $Lu=f$, for $f\in \Oka(U)$, has solutions $u\in
\Oka(U')$.

The problem that we shall consider in this paper consists of finding
large classes of irreducible algebraic hypersurfaces (i.e.
 codimension one algebraic
subvarieties of $\bC^n$)
$$  \Gamma_1,\ldots,  \Gamma_p  \label{LabGamma} $$
containing $0\in \bC^n$, and multiplicities $\mu_1,\ldots,\mu_p$
such that, for every domain $0\in U\subset \Omega$, there exists a
subdomain $0\in U'\subset U$ with the property that the boundary
value problem
\begin{equation}\label{mixedC0}
\left\{
\begin{aligned}
& Lu=f\\
& (D^\beta (u-g))|_{\Gamma_j}=0,\quad j=1,\ldots, p,\ 0\leq
|\beta|<\mu_j
\end{aligned}
\right.
\end{equation}
has a unique solution $u\in \Oka(U')$ for every $f,g\in \Oka(U)$. In
this case, we shall say that \eqref{mixedC0}, which we shall refer
to as a {\it mixed Cauchy problem for $L$ (at $0$)}, is {\it well
posed}. The classical Cauchy-Kowalevsky theorem corresponds to the
case of a hyperplane
$$\Gamma=\Gamma_1:=\{(z_1 ,...,z_n) \in \bC ^n\colon z_1 =0\}$$
(so $p=1 $ and the multiplicity is $k) $ and  $Q_k (D)= D^{(k,0,...,0)} , $
see e.g. \cite{Rauc97}, p. 15. By an analytic change of variables the
Cauchy-Kowalewsky theorem can be generalized to
the case of initial conditions on a hypersurface
$$\Gamma=\Gamma_1:=\{z\colon R(z)=0\}$$
 which is non-singular at $0$ (i.e.
the conormal vector $\zeta:=(\partial R/\partial z)(0)$ is not $0$),
and  non-characteristic with respect to $L$ at $0$ (i.e.
$Q_k(\zeta)\neq 0 )  ,$  see \cite{Rauc97}, p. 22. These Cauchy
problems are well posed. In this paper, we shall be interested in
the more difficult case of {\it singular} hypersurfaces, which is
not covered by the Cauchy-Kowalevsky theorem. Our methods of proof
depend on arguments using homogeneous power series in combination
with new decompositions of homogeneous polynomials, known as Fischer
decompositions; for more details we refer the reader to Section
\ref{Fischer}.

Before stating our main results and discussing previous results
along these lines, let us remark that it suffices to consider only
the mixed Cauchy problem with nulldata, i.e.
\begin{equation}\label{mixedC1}
\left\{
\begin{aligned}
& Lu=f\\
& (D^\beta u)|_{\Gamma_j}=0,\quad j=1,\ldots, p,\ 0\leq
|\beta|<\mu_j.
\end{aligned}
\right.
\end{equation}
since the equation $Lu=f$ has a solution in
$\Oka(U')$ for every $f\in \Oka(U)$.
For the remainder of
this paper, we shall consider only the problem \eqref{mixedC1} with
nulldata.

 In what follows, we shall give some equivalent reformulations of the
 mixed Cauchy
problem that will be more convenient from a technical point of view.
 Since each $\Gamma_j$ is an irreducible algebraic
hypersurface in $\bC^n$ with $0\in \Gamma_j$, there is an
irreducible polynomial $R_j(z)$, uniquely determined up to a
multiplicative constant and with $R_j(0)=0$, such that
$\Gamma_j:=\{z\colon R_j(z)=0\}$. The condition that
$$D^\beta u|_{\Gamma_j}=0,\quad 0\leq |\beta|<\mu_j,$$ for $u\in
\Oka(\Omega')$ is equivalent to $R_j^{\mu_j}$ dividing $u$ in the
ring $\Oka(\Omega')$, henceforth denoted
 by
$R_j^{\mu_j}|u$. Thus, if
we set $P:=R_1^{\mu_1}\ldots R_p^{\mu_p}$, then the mixed Cauchy
problem \eqref{mixedC1} can be equivalently formulated as follows,
\begin{equation}\label{mixedC2}
\left\{
\begin{aligned}
& Lu=f\\
& P| u.
\end{aligned}
\right.
\end{equation}
We shall refer to the polynomial $P$ in \eqref{mixedC2} as the {\it
divisor} in the mixed Cauchy problem. Now, given $f\in \Oka(U)$, a
function $u\in \Oka(U')$ is a solution to \eqref{mixedC1} or,
equivalently, to \eqref{mixedC2} if and only if $u=Pq$, for some
$q\in \Oka(U')$, and $L(Pq)=f$ in $U'$.
 In particular, {\it the mixed
Cauchy problem for the operator $L$ and divisor $P$ is well posed if
and only if, for every domain $0\in U\subset \Omega$, there is a
subdomain $0\in U'\subset U$ such that there exists a unique
solution $q\in\Oka(U')$ to the equation
\begin{equation}\label{mixedC3}
L(Pq)=f,\end{equation}
for every $f\in \Oka(U)$.} We shall use this formulation
of the mixed Cauchy problem in our main results,
Theorem
\ref{main0}, Theorem \ref{ThmMain2} and
Theorem \ref{abstract}.

 Previous results on the mixed Cauchy problem (at $0$)
includes a theorem by H\"ormander (\cite{Horm}, Theorem 9.4.2) in
the case where the divisor  $P(z)$ is a monomial $z^\gamma$ of degree
$|\gamma|=k$, $Q_k(D)=D^\gamma$, and certain sufficiently small
perturbations are allowed even in the principal part of $L$ (i.e.
$k_0$ in \eqref{defLL} is allowed to be $k$, but the coefficient
$a_\gamma(z)$ must be identically $0$ and there is a ``smallness"
requirement for those coefficients $a_\alpha(z)$ for which
$|\alpha|=k$). An early version of this theorem in two dimensions
goes back to Goursat (see e.g. \cite{Shap89}). Another, more recent
result is due to the first author, jointly with H. S. Shapiro
(\cite{EbSh96}, Theorem 3.1.1): There exists a number $k_0<k$
depending on $Q_k$
such that the mixed Cauchy problem with divisor
$P(z)=Q^*_k(z)=\overline{Q_k(\bar z)}$ has a unique solution $u\in
U'$ for every $f\in \Oka(U)$.

In this paper, we shall prove a result
(see Theorem \ref{main0}) for the mixed Cauchy problem for differential
operators $L$ of the type
 $Q_k(\zeta)=(B(\zeta))^m$ where  $k=2m$ and $B(\zeta)$ is a
nondegenerate quadratic form.  For the divisor $P(z)$ we only require that
it is a homogeneous polynomial
of degree $2m$ that is $B(\zeta)$-elliptic (see below for the
definition).
This result does not contain, nor is it contained in
the results from \cite{EbSh96} mentioned above. The results in
\cite{EbSh96} allow
 a more general
class of principal symbols $Q_k(\zeta)$, but, on the other hand, for
each $Q_k(\zeta)$ there is only one divisor $P(z)$ that can be used
in the mixed Cauchy problem, namely $Q^*_k(z)$. The result in the
present paper treats a smaller class of principal symbols, but for
each such principal symbol $Q_k(\zeta)$ there is a large class of
$P(z)$ that may be used as a divisor. We also give a more precise
result in $\bR^n$ for operators with the iterated Laplacian as their
principal symbol. The additional precision in this theorem concerns
the relation between $U'$ and $U$ (see Theorem \ref{ThmMain2}).

The paper is organized as follows. The main results are stated in
Section \ref{mainresults} and it is explained how Theorem
\ref{main0} follows from Theorem \ref{ThmMain2}. Section
\ref{examples} discusses two examples as illustrations of our main
results. The section after that introduces an integral that will be
used throughout the paper. The Fischer norms are then introduced and
some basic estimates are proved in Section \ref{Fischer}. The next
section contains further estimates and, in particular, the key
estimate (Theorem \ref{ThmHomA}) needed to prove Theorem
\ref{ThmMain2}. In the last section, Section \ref{proof}, we state
and prove a general result about mixed Cauchy problems in $\bR^n$
(Theorem \ref{abstract}), which together with the estimate in
Theorem \ref{ThmHomA} proves Theorem \ref{ThmMain2}.

\section{Main results}\label{mainresults}

In order to state our first main result, we need the following
definition.  Let $B(\zeta)$ be a nondegenerate quadratic form in
$\bC^n$, i.e.\ $B(\zeta)=\zeta^tB\zeta$ for some invertible,
symmetric $n\times n$ matrix with complex coefficients. By standard
linear algebra, there exists an invertible $n\times n$ matrix $A$
such that $B(A\tau)$ is equal to the standard nondegenerate
quadratic form $\Sigma(\tau)$,
\begin{equation}\label{standard}
\Sigma(\tau ):=\sum_{j=1}^n\tau_j^2.
\end{equation}
Let now $P_{2p}(z)$ be a homogeneous polynomial of degree $k:=2p\geq
2$. We shall say that $P_{2p}$ is {\it $B$-elliptic} if, for {\it
some} invertible $n\times n$ matrix $A$ such that
$B(A\tau)=\Sigma(\tau)$, the polynomial $P_{2p}(A^{-t}x)$ is
real-valued for $x\in \bR^n$ and there is a constant $\delta>0$ such
that
\begin{equation}\label{B-elliptic}
P_{2p}(A^{-t}x)\geq \delta (B(Ax))^p=\delta|x|^{2p},\quad x\in
\bR^n.
\end{equation}
Here $A^{-t} $ is the transpose of the inverse matrix $A^{-1} . $
For instance, if $B(\zeta)=\Sigma(\zeta)$  and $P_{2p}(x)$ is {\it
elliptic} in the usual sense, i.e.\ $P_{2p}(x)$ is real and
satisfies $P_{2p}(x)\geq\delta|x|^{2p}$ for $x\in \bR^n$, then of
course $P_{2p}$ is $B$-elliptic. However, we point out that $P_{2p}$
can be $\Sigma$-elliptic, even if $P_{2p}(x)$ fails to be elliptic,
as is illustrated by the following example.

\begin{example}
{\rm Let $\xi\in \bR$ and consider the following homogeneous
polynomial of degree $4$,
\begin{multline}
P(z)=P_4(z):=(\xi^4+(1+\xi^2)^2)z_1^4+(\xi^4+(1+\xi^2)^2)z_2^4
-12\xi^2(1 +\xi^2)z_1^2z_2^2\\
+4i\xi\sqrt{1+\xi^2}(1+2\xi^2)(z_1z_2^3+z_1^3z_2).
\end{multline}
The polynomial $P(x)$ is not real for $x\in \bR^n$ and, hence, is
not elliptic (nor is its real part elliptic if, say, $|\xi|\geq 1$).
However, if we let $A$ be the matrix
\begin{equation}
A:=\left(
\begin{matrix} i\xi & -\sqrt{1+\xi^2}\\\sqrt{1+\xi^2} & i\xi
\end{matrix}\right)
\end{equation}
then one can check that $\Sigma(A\tau)=\Sigma(\tau)$ and
$P(A^{-t}x)=x_1^4+x_2^4$. Since $P(A^{-t}x)$ is real and satisfies
$P(A^{-t}x)\geq \delta|x|^4$, we conclude that $P$ is
$\Sigma$-elliptic. }
\end{example}

We also mention that a homogeneous polynomial $P_{2p}(z)$ of degree
$2p$ is $B$-elliptic, for a given nondegenerate quadratic form
$B(\zeta)$, if and only if there exists a linear change of
coordinates $z=A^{-t} w$ such that $Q_{2p}(\partial/\partial
z):=(B(\partial/\partial z))^p$ in the new coordinates $w$ becomes
$\widetilde Q_{2p}(\partial/\partial w)=\Delta_\bC^p$, where
\begin{equation}\label{Qstandard}
\Delta_\bC:=\sum_{j=1}^n\frac{\partial^2}{\partial w_j^2},
\end{equation}
and the polynomial $\tilde P_{2p}(w):=P_{2p}(A^{-t}w)$ is elliptic
in the usual sense.

Our first main result is the following.

\begin{theorem}\label{main0} Let $B(\zeta)$ be a nondegenerate quadratic
form in $\bC^n$ and $p$ an integer $\geq 1$. Let $k:=2p$,
$Q_k(\zeta):=(B(\zeta))^p$, and consider the holomorphic partial
differential operator $L$ given by \eqref{defLL} with $k_0=p=k/2$.
Suppose $P(z)=P_k(z)$ is a homogeneous polynomial of degree $k=2p$
that is $B(\zeta)$-elliptic. Then, for any domain $0\in
U\subset\Omega$, there is a subdomain $0\in U'\subset U$ such that
the mixed Cauchy problem $$ L(Pq)=f$$ has a unique solution $q\in
\Oka(U')$ for every $f\in \Oka(U)$.
\end{theorem}

In the setting of Theorem \ref{main0}, as we mentioned above, we may
assume, possibly after a linear change of coordinates, that
$Q_{2p}(D)=\Delta_\bC^p$ and $P_{2p}(x)\geq \delta|x|^{2p}$ for
$x\in \bR^n$. Let $\Delta$ denote the usual Laplace operator in
$\bR^n$,
\begin{equation}\label{lap}
\Delta:=\sum_{j=1}^n\frac{\partial^2}{\partial x_j^2}.
\end{equation}
Theorem \ref{main0} will follow from a result about a mixed Cauchy
type problem in $\bR^n$ for partial differential operators whose
principal symbol is the iterated Laplace operator $\Delta^p$. To
formulate this result, we must introduce some more notation. Let
$B_{R}:=\left\{ x\in \mathbb{R}^{n}:\left| x\right| <R\right\} $ be
the open unit ball in $\mathbb{R}^{n}$ (where $0<R\leq \infty ).$ We
consider the algebra $A\left( B_{R}\right) $ of all infinitely
differentiable
functions $f:B_{R}\rightarrow \mathbb{C}$ such that for any compact subset
 $
K\subset B_{R}$ the homogeneous Taylor series $\sum_{m=0}^{\infty
}f_{m}\left( x\right) $ converges absolutely and uniformly to $f$ on
$K$; here, $f_{m}$ is the homogeneous polynomial of degree $m$
defined by the Taylor series of $f$
\begin{equation*}
f_{m}\left( x\right) =\sum_{\left| \alpha \right| =m}\frac{1}{\alpha
 !}\frac{%
\partial ^{\alpha }f}{\partial x^{\alpha }}\left( 0\right) x^{\alpha }.
\end{equation*}
Note that the functions in $A(B_R)$ are real-analytic. In fact, it
is known that $A\left( B_{R}\right) $ is isomorphic to
$\Oka(\widehat B_R)$, where $\widehat B_R\subset \bC^n$ denotes the
Lie ball of radius $R$
\begin{equation*}
\widehat{B_{R}}:=\left\{ z\in \mathbb{C}^{n}:\left| z\right| ^{2}+\sqrt{
\left| z\right| ^{4}-\left| z_1^2+\ldots + z_n^2\right|
^{2}}<R^{2}\right\},
\end{equation*}
and the isomorphism $\phi\colon \Oka(\widehat {B_R})\to A\left(
B_{R}\right)$ is simply given by $\phi(f):= f|_{B_R}$. (See
\cite{Sici74} for this result; see also \cite{Rend05}, Section 8.)
We observe that the isomorphism $\phi$ commutes with differentiation
in the following way
$$
\phi\left(\frac{\partial^{|\alpha|}f}{\partial z^\alpha}\right)=
\frac{\partial^{|\alpha|}\phi(f)}{\partial x^\alpha}.
$$
Since any domain $0\in U$ contains a Lie ball of some radius and
every Lie ball contains an open neighborhood $U'$ of $0$, we
conclude, as claimed above, that Theorem \ref{main0} indeed is a
consequence of the following result in $\bR^n$.

\begin{theorem}
\label{ThmMain2} Suppose that $P_{2p}\left( x\right) $ is
homogeneous of degree $2p$ and elliptic, i.e.\ there exists
$\delta>0$ such that $ P_{2p}\left( x\right) \geq \delta\left|
x\right| ^{2p} $ for all $x\in \mathbb{R}^{n}.$ Let $0\leq k_{0}\leq
p$ be an integer, $R>0$ a positive number, $a_\alpha(x)$ functions
in $A(B_{R})$ for every multi-index $\alpha\in \bN_0^n$ with
$|\alpha|\leq k_0$, and
\begin{equation*}
L=\Delta ^{p}+\sum_{|\alpha|\leq k_0}a_\alpha(x) D^\alpha.
\end{equation*}
If $k_{0}<p$ then the operator $q\mapsto L\left(P_{2p}\, q\right) $
is a bijection from $A\left( B_{R}\right) $ onto $A\left(
B_{R}\right)$. If $k_{0}= k$ then there exists $r>0$ such that the
equation $ L\left(P_{2p}\,q\right)=f $ has a unique solution $q\in
A\left( B_{r }\right) $ for every  $f\in A\left( B_{R}\right)$.
\end{theorem}

The proof of Theorem \ref{ThmMain2} hinges on new estimates for a
real version of the Fischer norm (see Theorem \ref{ThmHomA}) that go
back to the paper \cite{Rend05} by the second author. Theorem
\ref{ThmMain2} follows then from a general result (Theorem
\ref{abstract}) about real mixed Cauchy type problems. The latter
theorem is analogous to a similar theorem about complex Cauchy
problems in \cite{EbSh96}.

We note that if $Q_{2p}(D)=\Delta_\bC^p$ in Theorem \ref{main0} (as
we may assume), then the homogeneous polynomial
$P_{2p}(x)=Q_{2p}^*(x)=|x|^{2p}$ is $B$-elliptic. Thus, both Theorem
\ref{main0} and Theorem 3.1.1 in \cite{EbSh96} apply to the mixed
Cauchy problem for $L$ given by \eqref{defLL} with divisor
$P_{2p}(z)=\sum z_j^2$. In this particular situation, the result in
\cite{EbSh96} is more general: the number $k_0$ in \eqref{defLL} can
be chosen to be $3p/2$ (see \cite{EbSh96}, p. 261), whereas in the
present paper only $k_0=p$ is allowed. The reason for this is that
\cite{EbSh96} utilizes the complex Fischer norm, rather than the
real one used in this paper, and when $Q_{2p}(D)=\Delta_\bC^p$,
$P_{2p}(z)=\sum z_j^2$, a stronger estimate holds for the complex
Fischer norm (see Subsection \ref{comments}). The advantage of the
real norm, of course, is that it allows a much more general class of
divisors.

\section{Examples and applications}\label{examples}

In this section, we apply Theorem \ref{ThmMain2} to a couple of
explicit examples.
Before proceeding, we should perhaps point out that, in general, the
mixed Cauchy problem for $L$ with divisor $P$ is {\it not} well
posed, even if $P$ is a homogeneous polynomial of degree $k$, as is
illustrated by the following simple example.

\begin{example}\label{nosol} {\rm Consider the
 complex ``Laplace operator"
 in
two variables
$$
L=\Delta_\bC:=\frac{\partial^2 }{\partial z_1^2}+\frac{\partial^2
}{\partial z_2^2}
$$
and the homogeneous polynomial $P(z)=z_1z_2$. Note that $q=1$ solves
$L(Pq)=0$ and, hence, the uniqueness fails. It is also easy to see
that the equation $L(Pq)=1$ has no solution in any neighborhood of
the origin. }
\end{example}

In \cite[Section 5]{EbSh96}, it is also shown that solvability can
fail even when uniqueness holds. For instance, if we take $L$ to be
the complex Laplace operator in $\bC^2$ and
$$P(z)=z_2(z_1^2+(z_2-1)^2-1)=z_2(z_1^2+z_2^2-2z_2),$$
 then uniqueness holds in the mixed
Cauchy problem at $0$ but $L(P q)=f$ is in general not solvable.

 The
problem of deciding for which polynomials $P(z)$ the mixed Cauchy
pro\-blem $L(P q)=0$, where $L$ is the com\-plex Laplace operator, has
$q=0$ as its unique solution has been addressed in e.g.\ \cite{A},
\cite{AK}.

\begin{example} \label{ex1} {\rm Consider the holomorphic
partial
 differential operator
\begin{equation}\label{Lex1}
L:=\Delta_\bC^2+\sum_{j=1}^n a_n(z)\frac{\partial}{\partial
z_j}+b(z),
\end{equation}
where, for simplicity, the coefficients $a_j(z)$ and $b(z)$ are
assumed to be entire functions in $\bC^n$. Let $P(z)=\sum_{j=1}^n
z_j^4$ and note that $P(x)\geq \delta|x|^4$ for $x\in \bR^n$. Since
$k_0=1<2$, it follows from Theorem \ref{ThmMain2} and the remarks
preceding it that the mixed Cauchy problem
$$
L(Pq)=f
$$
has a unique solution $q\in \Oka(\widehat{B_R})$ for any $f\in
\Oka(\widehat{B_R})$ and any $R>0$. This illustrates the fact that
if $U$ is a Lie ball, then one can take $U'=U$ in Theorem
\ref{main0} provided that $k_0<p$ (and the coefficients are analytic
in $\widehat {B_{2R}}$). }
\end{example}

\begin{example} {\rm Let $\square$ denote the wave operator in
$\bR^n\times \bR$, $$ \square:=\sum_{j=1}^n
\frac{\partial^2}{\partial x_j^2}-\frac{\partial^2}{\partial t^2}$$
and consider the real partial differential operator
\begin{equation}
L:=\square+a(x,t),
\end{equation}
where $a(x,t)$ is, say, in $A(\bR^{n+1})$. Let
$$
P(x,t):=\sum_{j=1}^nx_j^2-t^2,
$$
so that $\{(x,t)\colon P(x,t)=0\}$ is the light cone. Observe that
the linear change of variables $y=it$ transforms $\square$ into the
Laplace operator $\Delta$ in $\bR^{n+1}$ and $\widetilde
P(x,y):=P(x,it)$ becomes
$$
\widetilde P(x,y)=\sum_{j=1}^nx_j^2+y^2,
$$ which is clearly elliptic. An application of Theorem \ref{main0}
and the remark made in Example \ref{ex1} above (here, $k_0=0<1=p$)
yields (the probably well known result) that the real Cauchy problem
\begin{equation}\label{wave}
L(P q)=f
\end{equation}
has a unique solution $q$ in $A(D_R)$ for every $f\in A(D_R)$. Here,
$D_R$ is the real domain
\begin{equation*}
D_{R}:=\left\{ (x,y)\in \mathbb{R}^{n}\times \bR:\left| x\right| ^{2}
 +|y|^2+\sqrt{
\left( |x^2|+|y|^2\right) ^{2}-\left| x_1^2+\ldots x_n^2-y^2\right|
^{2}}<R^{2}\right\},
\end{equation*}
and $A(D_R)$ denotes the restriction to $D_R$ of functions that are
holomorphic in
$$
\left\{ (z,w))\in \mathbb{C}^{n}\times \bC:\left| z\right| ^{2}=
 +|w|^2+\sqrt{%
\left( |z^2|+|w|^2\right) ^{2}-\left| z_1^2+\ldots z_n^2-w^2\right|
^{2}}<R^{2}\right\}
$$
We point out that the light cone, which carries the null data in
\eqref{wave}, is everywhere characteristic for the wave operator
$\square$.
 }
\end{example}

\section{A special integral}

Throughout the paper we shall use frequently the following notation:
\begin{equation*}
I_{m}:=\int_{0}^{\infty }e^{-r^{2}}r^{m}dr\text{ for }m\in \mathbb{N}_{0}.
\end{equation*}
This integral is well known, and for the even case (see p. 265 in
 \cite{Rem}
) we have
\begin{equation}
I_{2m}=\frac{\sqrt{\pi }}{2}\frac{\left( 2m\right)
 !}{m!}2^{-2m}=\frac{\sqrt{%
\pi }}{2}\frac{1\cdot 3\cdot 5\cdot ...\cdot \left( 2m-1\right) }{2^{m}}\leq
m!,  \label{deI2m}
\end{equation}
while in the odd case a simple substitution argument gives
\begin{equation}
I_{2m+1}=\int_{0}^{\infty
 }e^{-r^{2}}r^{2m+1}dr=\frac{1}{2}\int_{0}^{\infty
}e^{-x}x^{m}dx=\frac{1}{2}m!.  \label{eqI2m1}
\end{equation}
We shall use the following identity.

\begin{proposition}\label{Iest}
For positive integers $m,k,j,n$,
\begin{equation}
\frac{I_{2m+2jk+n-1}}{I_{2m+n-1}}=\frac{1}{2^{jk}}\left(
n+2m\right)(n+2m+2)\ldots\left(n+2m+2jk-2\right).
\label{EEE1}
\end{equation}
\end{proposition}

\begin{proof}
First assume that $n-1$ is even and write $n-1=2l.$ Then by
\eqref{deI2m}
\begin{eqnarray*}
\frac{I_{2m+2jk+n-1}}{I_{2m+n-1}} &=&\frac{1}{2^{jk}}\frac{1\cdot
3\cdot 5\cdot ...\cdot \left( 2\left( m+jk+l\right) -1\right)
}{1\cdot 3\cdot 5\cdot
...\cdot \left( 2\left( m+l\right) -1\right) } \\
&=&\frac{1}{2^{jk}}\left( 2m+2l+1\right)\left( 2m+2l+3\right) \ldots\left
( 2m+2l+2jk-1\right) \\
&=&\frac{1}{2^{jk}}\left( 2m+n\right)(2m+n+2) \ldots \left(
2m+n-2+2jk\right) .
\end{eqnarray*}
If $n-1$ is odd, then write $n=2l.$ We obtain
\begin{equation}
\frac{I_{2m+2jk+n-1}}{I_{2m+n-1}} =\frac{I_{2m+2jk+2l-1}}{I_{2m+2l-1}}=\frac{
\left( m+jk+l-1\right) !}{\left( m+l-1\right) !}.
\label{EEE2}
\end{equation}
On the other hand, the right hand side of (\ref{EEE1}) for  $n=2l $ is equal to
\begin{equation*}
  \frac{\left( 2l+2m\right)(2l+2m+2)\ldots
 \left( 2l+2m+2jk-2\right)}{2^{jk}}
\end{equation*}
which is equal to
$
\left( l+m\right)(l+m+1)\ldots\left(l+m+ jk-1\right)
$. In view of (\ref{EEE2}) the proof is finished.
\end{proof}

\section{Basic estimates in Fischer type spaces}\label{Fischer}

Let $\mathbb{C}\left[ x_{1},...,x_{n}\right] $ be the space of all
polynomials in $n$ variables with complex coefficients. An important inner
product on $\mathbb{C}\left[ x_{1},...,x_{n}\right] $ is the so-called
 \emph{
Fischer inner product}, or the \emph{apolar inner product}, defined by
\begin{equation*}
\left\langle P,Q\right\rangle _{F}:=\sum_{\alpha \in \mathbb{N}
_{0}^{n}}\alpha !c_{\alpha }\overline{d_{\alpha }}
\end{equation*}
for polynomials $P\left( x\right) =\sum_{\left| \alpha \right| \leq
N}c_{\alpha }x^{\alpha }$ and $Q\left( x\right) =\sum_{\left| \alpha
\right| \leq N}d_{\alpha }x^{\alpha }$, which has been used by
several authors, see e.g. in chronological order
 \cite{Fisc17}, \cite{Barg61}, \cite{CaZy64},
\cite{Kura66}, \cite{NeSh66}, \cite{NeSh68}, \cite{Dono},
\cite{Kura71}, \cite{StWe71}, \cite{Shap89}, \cite {BBEM90},
\cite{deRo92}, \cite{EhRo93}(and the references given there),
\cite{Zeil}, \cite{EbSh95}, \cite{BeDe95}, \cite {Beau97},
\cite{Thar97}, \cite{Vegt99}, \cite{FGS98}, \cite{Koun00}, and \cite
{ArGa01}.
 This inner product has
the property that the adjoint map of the
differentiation operator $Q\left( D\right) $ is the multiplication operator
 $
M_{Q^{\ast }}$, defined by $M_{Q^{\ast }}\left( f\right) =Q^{\ast }\cdot
 f$;
so this means that
\begin{equation}
\left\langle Q\left( D\right) f,g\right\rangle _{F}=\left\langle
f,Q^{\ast }\cdot g\right\rangle _{F}=\left\langle f,M_{Q^{\ast
}}g\right\rangle _{F} \label{eqFischeradjoint}
\end{equation}
for all polynomials $f,g\in \mathbb{C}\left[ x_{1},...,x_{n}\right] $ where
 $
Q^{\ast }$ is the polynomial obtained by conjugating the coefficients. It
was already observed by V. Barg\-mann in 1961 (see \cite{Barg61}) that
\begin{equation}
\left\langle f,g\right\rangle _{F}=\frac{1}{\pi ^{n}}\int_{\mathbb{R}
^{n}}\int_{\mathbb{R}^{n}}f\left( x+iy\right) \overline{g\left( x+iy\right)}
e^{-\left| x\right| ^{2}-\left| y\right| ^{2}}dxdy=\frac{1}{\pi
^{n}}\int_{\mathbb{C}^{n}} f\left( z\right)
\overline{g(z)}e^{-\left| z\right| ^{2}}dA_z \label{Bargmann}
\end{equation}
where $dx$, $dy$ denote the Lebesgue measure on $\mathbb{R}^{n}$ and
$dA_z$ the Lebesgue measure on $\bC^n\cong\bR^{2n}$. In passing, we
note that the \emph{Bargmann space} $\mathcal{F}_{n}$ (also called
\emph{Fock} or
\emph{Fischer space}) is defined as the space of all entire functions $f:
\mathbb{C}^{n}\rightarrow \mathbb{C}$ which satisfy
\begin{equation*}
\left\| f\right\| _{F}^{2}=\frac{1}{\pi
^{n}}\int_{\mathbb{C}^{n}}\left| f\left( z\right) \right|
^{2}e^{-\left| z\right| ^{2}}dA_z<\infty.
\end{equation*}

In analogy with equation (\ref{Bargmann}), we shall consider the
following real version of the Fischer inner product:
\begin{equation}
\left\langle f,g\right\rangle _{rF}:=\int_{\mathbb{R}^{n}}f\left( x\right)
\overline{g\left( x\right) }e^{-\left| x\right| ^{2}}dx,
\label{realBargmann}
\end{equation}
which has been useful for solving the Hayman conjecture for uniqueness sets
of polyharmonic functions and for solving the Khavinson-Shapiro conjecture
for the Dirichlet problem, see \cite{Rend05}. Note that in (\ref
{realBargmann}) we consider a polynomial as a function on the space
 $\mathbb{
R}^{n}$, while in (\ref{Bargmann}) it is considered as a function on the
space $\mathbb{C}^{n}.$

We should point out that the two inner products have some important
differences, e.g. the adjoint map for the multiplication operator
$M_{Q}$ for the inner product $\left\langle \cdot ,\cdot
\right\rangle _{rF}$ is not the differentiation operator but just
the operator $M_{Q^{\ast }}.$ However, it is a somewhat surprising
fact that the two inner product share many properties as well. As an
illustrative example we begin with the following proposition, part
of which will be crucial in the proof of Theorem \ref{abstract}
below.

\begin{proposition}\label{Pest} Let $k$ and $n$ be positive integers.
Let $\bS^{n-1}$ denote the unit sphere in $\bR^n$ and
$\Sigma^{2n-1}$ the unit sphere in $\bC^n\cong\bR^{2n}$. Let $P_{k}$
be a homogeneous polynomial of degree $k$ in $n$ variables,
$M_{\bR}:=\max_{\theta\in \mathbb S^{n-1}}|P_k(\theta)|$, and
$M_{\bC}:=\max_{\eta\in \Sigma^{2n-1}}|P_k(\eta)|$. Then, for any
homogeneous polynomial $f_{m}$ of degree $m$ in $n$ variables,
\begin{equation}
\left\| P_{k}f_{m}\right\| _{F}\leq
M_{\bC}\sqrt{\frac{I_{2m+2k+2n-1}}{I_{2m+2n-1}}} \left\|
f_{m}\right\| _{F},\ \left\| P_{k}f_{m}\right\| _{rF}\leq
M_{\bR}\sqrt{\frac{I_{2m+2k+n-1}}{I_{2m+n-1}}}\left\| f_{m}\right\|
_{rF} \label{eqmult}.
\end{equation}
In particular, for fixed $k$ and $n$, there are constants
$C_{k,n}>0$ and $D_{k,n}>0$ such that
\begin{equation}
\left\| P_{k}f_{m}\right\| _{F}\leq  C_{k,n}M_{\bC} \sqrt{ (1+m)
^{k}}\left\| f_{m}\right\| _{F},\ \left\| P_{k}f_{m}\right\|
_{rF}\leq  D_{k,n}M_{\bR} \sqrt{ (1+m) ^{k}}\left\| f_{m}\right\|
_{rF} \label{eqmult2}.
\end{equation}
\end{proposition}

\begin{proof}
Let us consider first the norm $\left\| {\cdot}\right\| _{rF}.$ By
introducing polar coordinates, it is easy to see that for a
homogeneous polynomial of degree $m$
\begin{equation*}
\left\| f_{m}\right\| _{rF}^{2}=\left\langle f_m,f_m\right\rangle
_{rF}=I_{2m+n-1}\int_{\mathbb{S}^{n-1}}\left| f_m\left( \theta
\right) \right| ^{2}d\theta .
\end{equation*}
Applied to $P_{k}f_{m}$ this gives
\begin{equation*}
\left\| P_{k}f_{m}\right\| _{rF}^{2}=I_{2m+2k+n-1}\int_{\mathbb{S}
^{n-1}}\left| P_{k}\left( \theta \right) f_m\left( \theta \right)
\right| ^{2}d\theta .
\end{equation*}
 Then
\begin{equation*}
\left\| P_{k}f_{m}\right\| _{rF}^{2}\leq \frac{I_{2m+2k+n-1}}{I_{2m+n-1}}
M_{\bR}^{2}\left\| f_{m}\right\| _{rF}^{2}.
\end{equation*}
This proves the second inequality in (\ref{eqmult}). From
Proposition \ref{Iest}, it is easy to see that
\begin{equation*}
\frac{I_{2m+2k+n-1}}{I_{2m+n-1}}\leq (m+k+n/2-1) ^{k}.
\end{equation*}
Clearly, for fixed $k$ and $n$, there exists a constant $D_{k,n}$
such that the second inequality in \eqref{eqmult2} holds.

For the computation of the norm $\left\| {\cdot}\right\| _{F}$, we
note that
\begin{equation*}
\left\| f_{m}\right\| _{F}^{2}=I_{2m+2n-1}\int_{\Sigma
^{2n-1}}\left| f_m\left( \eta \right) \right| ^{2}d\eta .
\end{equation*}
 Then
\begin{equation*}
\left\| P_{k}f_{m}\right\| _{F}^{2}\leq \frac{I_{2m+2k+2n-1}}{I_{2m+2n-1}}
{M_{\bC}}^{2}\left\| f_{m}\right\| _{F}^{2}.
\end{equation*}
This proves the first inequality in \eqref{eqmult}.  The first
inequality in \eqref{eqmult2} follows easily from Proposition
\ref{Iest} as above.
\end{proof}

As a second example, also used in the proof of Theorem
\ref{abstract}, we consider estimates of the derivative of
homogeneous polynomials:

\begin{proposition}
\label{PropDeriv}Let $\alpha \in \mathbb{N}_{0}^{d}$ be a multi-index and $
D^{\alpha }$ be the corresponding differential operator. Then
\begin{equation*}
\left\| D^{\alpha }f_{m}\right\| _{F}\leq \sqrt{m^{\left| \alpha \right| }}
\left\| f_{m}\right\| _{F}\text{ and }\left\| D^{\alpha }f_{m}\right\|
_{rF}\leq \sqrt{\left( 2m\right) ^{\left| \alpha \right| }}\left\|
f_{m}\right\| _{rF}
\end{equation*}
for any homogeneous polynomial $f_{m}$ of degree $m$.
\end{proposition}

\begin{proof}
By a simple induction argument, it is sufficient to prove the
statement for the differential operator $D_{j}:=\frac{\partial
}{\partial x_{j}}.$ In case of $\left\| {\cdot}\right\| _{F}$ we
repeat (for convenience of the reader) the argument already given in
\cite{KhSh92} (or see \cite [p. 256]{EbSh95}): By Euler's formula
one has
\begin{equation*}
\sum_{j=1}^{n}z_{j}D_{j}f_{m}=mf_{m}.
\end{equation*}
Taking the Fischer inner product with $f_{m},$ and using that multiplication
by $z_{j}$ is adjoint to $D_{j}$ one obtains
\begin{equation*}
\sum_{j=1}^{n}\left\| D_{j}f_{m}\right\| _{F}^{2}=m\left\| f_{m}\right\|
_{F}^{2}.
\end{equation*}
In particular,
\begin{equation*}
\left\| D_{j}f_{m}\right\| _{F}\leq \sqrt{m}\left\| f_{m}\right\| _{F}.
\end{equation*}

Note that the previous argument does not apply to the norm $\left\|
\cdot \right\| _{rF}$ since $D_{j}$ is not the adjoint of $z_{j}.$
However, a
simple argument using partial integration shows that for $j=1,...,n$ and
 $
f,g\in \mathbb{C}\left[ x_{1},...,x_{n}\right] $
\begin{equation*}
\left\langle \frac{\partial }{\partial x_{j}}f,g\right\rangle
_{rF}+\left\langle f,\frac{\partial }{\partial x_{j}}g\right\rangle
_{rF}=2\left\langle x_{j}\cdot f,g\right\rangle _{rF}.
\end{equation*}
Replace $f$ by $\frac{\partial }{\partial x_{j}}f$, and sum up, then
\begin{equation*}
\left\langle \Delta f,g\right\rangle _{rF}+\sum_{j=1}^{n}\left\langle
 \frac{%
\partial }{\partial x_{j}}f,\frac{\partial }{\partial x_{j}}g\right\rangle
_{rF}=2\sum_{j=1}^{n}\left\langle x_{j}\frac{\partial }{\partial x_{j}}
f,g\right\rangle _{rF}.
\end{equation*}
For a homogeneous polynomial $f_m$ of degree $m$ Euler's formula
yields
\begin{equation}
\left\langle \Delta f_m,f_m\right\rangle _{rF}+\sum_{j=1}^{n}\left\|
 \frac{%
\partial }{\partial x_{j}}f_m\right\| _{rF}^{2}=2m\left\langle
f_m,f_m\right\rangle _{rF}.  \label{eqnabla1}
\end{equation}
Hence it suffices to show that $\left\langle \Delta
f_m,f_m\right\rangle _{rF}\geq 0$, which will be done in the next
proposition.
\end{proof}

\begin{proposition}\label{Prop10}
For any homogeneous polynomial $f_m$ of degree $m$
\begin{equation*}
\left\langle \Delta f_m,f_m\right\rangle _{F}=0\text{ and
}\left\langle \Delta f_m,f_m\right\rangle _{rF}\geq 0
\end{equation*}
\end{proposition}

\begin{proof}
The identity for $\langle\cdot,\cdot\rangle_F$ is trivial since
$\Delta f$ is polynomial of degree $m-2$ and homogeneous polynomials
of different degree are always orthogonal for the Fischer inner
product. It is an elementary fact that, for a homogeneous harmonic
polynomial $h\left( x\right) $, the following formula holds
\begin{equation}
\Delta \left( \left| x\right| ^{2s}h\left( x\right) \right) =2s\left[
2s-2+2\deg h+n\right] \cdot \left| x\right| ^{2s-2}h\left( x\right) .
\label{eqDelts}
\end{equation}
For the inequality for the real inner product
 $\langle\cdot,\cdot\rangle_{rF}$, we
consider the  {\it Gau\ss{} decomposition} of $
f_m$: there exist homogeneous harmonic polynomials $h_{m-2s}$ of
degree $m-2s$ such that $f_m=\sum_{s=0}^{N}\left| x\right|
^{2s}h_{m-2s}$ with $N=[m/2]$, see e.g. \cite{ABR92}, p. 76. Then,
according to (\ref{eqDelts}), $
\Delta \left( \left| x\right| ^{2s}h_{m-2s}\right) =c_{s}\left|
x\right| ^{2s-2}h_{m-2s},$ with
\begin{equation*}
c_{s}:=2s\left( 2s-2+n+2\deg h_{m-2s}\right) \geq 0.
\end{equation*}
Thus
\begin{equation*}
\left\langle \Delta f_m,f_m\right\rangle
=\sum_{s=0}^{N}\sum_{j=0}^{N}c_{s}\left\langle \left| x\right|
^{2s-2}h_{m-2s},\left| x\right| ^{2j}h_{m-2j}\right\rangle _{rF}.
\label{eqsumdelta2}
\end{equation*}
Furthermore,
\begin{equation*}
\left\langle \left| x\right| ^{2s-2}h_{m-2s},\left| x\right|
^{2j}h_{m-2j}\right\rangle _{rF}=I_{2m+n-3}\int_{\mathbb{S}%
^{n-1}}h_{m-2s}\left( \theta \right) h_{m-2j}\left( \theta \right)
d\theta .
\end{equation*}
Since $\deg h_{m-2s}-\deg h_{m-2j}=2\left( j-s\right) ,$ we see that
$\left\langle h_{m-2s},h_{m-2j}\right\rangle _{rF}=0$ for $s\neq j$.
Hence there is only a contribution in (\ref{eqsumdelta2}) for $s=j,$
and we obtain
\begin{equation*}
\left\langle \Delta f_m,f_m\right\rangle
_{rF}=\sum_{s=0}^{N}c_{s}\left\langle \left| x\right|
^{2s-2}h_{m-2s},\left| x\right| ^{2s}h_{m-2s}\right\rangle _{rF}\geq
0.
\end{equation*}
\end{proof}

\section{Operators acting on homogeneous polynomials}

Let $\mathcal{P}_{m}\left( \mathbb{R}^{n}\right) $ denote the space
of all homogeneous polynomials with complex coefficients of degree
$m$ in $n$
variables. We consider first the operator $F_{2p}:$ $\mathcal{P}
_{m}\left( \mathbb{R}^{n}\right) \rightarrow \mathcal{P}_{m}\left(=
 \mathbb{R}
^{n}\right) $ defined by
\begin{equation*}
F_{2p}\left( q\right) :=\Delta ^{p}\left( \left| x\right|
^{2p}q\right) .
\end{equation*}
A simple induction argument using the formula (\ref{eqDelts}) shows
that, for any homogeneous harmonic polynomial $h$,
\begin{equation}
F_{2p}\left( \left| x\right| ^{2s}h\right) =\Delta ^{p}\left( \left|
x\right| ^{2s+2p}h\right) =d_{p}\left( s,\deg h\right) \left|
x\right| ^{2s}h \label{eqeigen}
\end{equation}
where $d_{p}\left( s,m\right) $ is the number
\begin{equation}
d_{p}\left( s,m\right)= 2^{p}\left( s+p\right) ....\left( s+1\right)
\cdot \left( 2s+2p-2+n+2m\right) ....\left( 2s+n+2m\right) .
\label{eigen2}
\end{equation}
From this one obtains the following well-known result; for the
reader's convenience, we shall sketch the proof.

\begin{proposition}
The space $\mathcal{P}_{m}\left( \mathbb{R}^{n}\right) $ has a basis
consisting of eigenvectors for the operator
$F_{2p}:\mathcal{P}_{m}\left( \mathbb{R}^{n}\right) \rightarrow
\mathcal{P}_{m}\left( \mathbb{R}^{n}\right) $ such that the lowest
eigenvalue is greater
 than or equal to
\begin{equation}
e_{p,m}=2^{p}p!\left( 2m+n\right) \left( 2m+n+2\right) ...\left(
2m+n+2\left( p-1\right) \right) .  \label{eqmini2}
\end{equation}
\end{proposition}

\begin{proof}
Let $m\geq 1$ be fixed. Let $\mathcal{H}_{m-2s}\left( \mathbb{R}^{n}\right)
 $
be the space of all harmonic polynomials of degree $m-2s,$ and let $
Y_{m-2s,l}$ for $l=1,...,a_{m-2s}:=\dim \mathcal{H}_{m-2s}\left(
 \mathbb{R}
^{n}\right) $ be a basis of $\mathcal{H}_{m-2s}\left( \mathbb{R}^{n}\right)
. $ Then
\begin{equation}
\left| x\right| ^{2s}Y_{m-2s,l},\quad s=0,...,\left[ m/2\right] ,
l=1,...,a_{m-2s}  \label{eqbasis}
\end{equation}
are homogeneous polynomials of degree $m,$ and by (\ref{eqeigen})
they are clearly eigenfunctions of $F_{2p}$ with eigenvalue
$d_{p}\left( s,m-2s\right)$ and
\begin{equation*}
d_{p}\left( s,m-2s\right) :=2^{p}\left( s+p\right) ....\left(
s+1\right) \cdot \left( 2m-2s+2p-2+n\right) ....\left(
2m-2s+n\right) .
\end{equation*}
The minimal value for these numbers, ranging from $s=0,...,\left[
m/2\right] ,$ is attained for $s=0$ which gives (\ref{eqmini2}). The
Gau\ss{} decomposition of a polynomial
(see the proof of Proposition \ref{Prop10}) shows that (\ref{eqbasis}) is
indeed a basis of $\mathcal{P}_{m}\left( \mathbb{R}^{n}\right) $.
\end{proof}

\begin{proposition}\label{standest}
\label{PropDrF}For a homogeneous polynomial $f_{m}$ of degree $m$,
we have
\begin{equation*}
\left\| \Delta ^{p}\left( \left| x\right| ^{2p}f_{m}\right) \right\|
_{rF}\geq e_{p,m}\left\| f_{m}\right\| _{rF}.
\end{equation*}
\end{proposition}

\begin{proof}
Let $f_{m}=\sum_{s=0}^{N }\left| x\right| ^{2s}h_{m-2s}$, with
$N:=\left[ m/2\right]$, be the Gau\ss{} decomposition with harmonic
polynomials $h_{m-2s}$ of degree $m-2s$ for $s=0,...,N.$ We compute
the inner product $ \left\langle
F_{2p}f_{m},F_{2p}f_{m}\right\rangle _{rF}$ for $F_{2p}:=\Delta
^{p}\left( \left| x\right| ^{2p}\cdot \right) :$
\begin{equation*}
\left\langle F_{2p}f_{m},F_{2p}f_{m}\right\rangle
_{rF}=\sum_{s=0}^{N}\sum_{j=0}^{N}d_p(s,m-2s)
d_p(j,m-2j)\left\langle \left| x\right| ^{2s}h_{m-2s},\left|
x\right| ^{2j}h_{m-2j}\right\rangle _{rF}.
\end{equation*}
Since $\deg h_{m-2s}-\deg h_{m-2j}=2\left( j-s\right) ,$ we see that
$\left\langle h_{m-2s},h_{m-2j}\right\rangle _{rF}=0$ for $s\neq j$.
Hence
\begin{equation*}
\left\langle F_{2p}f_{m},F_{2p}f_{m}\right\rangle
_{rF}=\sum_{s=0}^{N}d_p(s,m-2s) ^{2}\left\langle \left| x\right|
^{2s}h_{m-2s},\left| x\right| ^{2s}h_{m-2s}\right\rangle _{rF}.
\end{equation*}
Similarly, $\left\langle f_{m},f_{m}\right\rangle
_{rF}=\sum_{s=0}^{N}\left\langle \left| x\right|
^{2s}h_{m-2s},\left| x\right| ^{2s}h_{m-2s}\right\rangle _{rF}.$
Hence
\begin{equation*}
\left\| F_{2p}f_{m}\right\| _{rF}\geq e_{p,m}\left\| f_{m}\right\|
_{F}.
\end{equation*}
\end{proof}

We shall now give the basic $\|\cdot\|_{rF}$-estimate for the
operator
\begin{equation*}
f\longmapsto \Delta ^{p}\left( P_{2p}\cdot f_{m}\right),
\end{equation*}
which will be used in the proof of Theorem \ref{ThmMain2}. We shall
show in the comments below (Subsection \ref{comments})
that the result is sharp even if $%
P_{2p}\left( x\right) =\left| x\right| ^{2p}.$ This is in contrast
with the case of the complex Fischer norm $\left\| {\cdot}\right\|
_{F}$, where a better estimate than (\ref{eqbas2}) holds for
$P_{2p}\left( z\right) =(\sum z_j^2) ^{p}$,
see (\ref{CFest}).

\begin{theorem}
\label{ThmHomA}Let $P_{2p}\left( x\right) $ be a homogeneous
polynomial of degree $2p$ and suppose that there is a $\delta>0$
such that $P_{2p}\left( x\right) \geq \delta\left| x\right| ^{2p}$
for all $x\in \mathbb{R}^{n}.$ Then there exists a constant $C_{1}$
such that, for each homogeneous polynomial $f_{m}$ of degree $m\in
\mathbb{N} _{0},$
\begin{equation}
\left\| \Delta ^{p}\left( P_{2p}\cdot f_{m}\right) \right\|
_{rF}\geq C_{1}e_{p,m}\left\| f_{m}\right\| _{rF}.  \label{eqbas1}
\end{equation}
Moreover, there exists a constant $C_{2}$ such that, for each
homogeneous polynomial $f_{m}$ of degree $m\in \mathbb{N}_{0}$,
\begin{equation}
\left\| \Delta ^{p}\left( P_{2p}\cdot f_{m}\right) \right\|
_{rF}\geq C_{2}\left\| P_{2p}\cdot f_{m}\right\| _{rF}.
\label{eqbas2}
\end{equation}
\end{theorem}

\begin{proof}
As above, we let  $F_{2p}\left( u\right) =\Delta ^{p}\left( \left|
x\right| ^{2p}u\right)$. Let $f_{m}$ be given and define
$g_{m}=:\Delta ^{p}\left( P_{2p}f_{m}\right) .$ Since $F_{2p}$ is a
bijection there exists $u_{m}$ such that $F_{2p}\left( u_{m}\right)
=\Delta ^{p}\left( \left| x\right| ^{2p}u_{m}\right) =g_{m}.$
Proposition \ref{PropDrF} yields
\begin{equation}
\left\| \Delta ^{p}\left( P_{2p}f_{m}\right) \right\| _{rF} =\left\|
g_{m}\right\| _{rF} =\left\| F_{2p}\left( u_{m}\right) \right\|
_{rF}\geq e_{p,m}\left\| u_{m}\right\| _{rF}.  \label{eqgg}
\end{equation}
Since obviously $\left\| \left| x\right| ^{2p}u_{m}\right\|
 _{rF}^{2}=\frac{
I_{2m+4p+n-1}}{I_{2m+n-1}}\left\| u_{m}\right\| _{rF}^{2}$ one
obtains
\begin{equation}
\left\| \Delta ^{p}\left( P_{2p}f_{m}\right) \right\| _{rF}\geq
 e_{p,m}\sqrt{%
\frac{I_{2m+n-1}}{I_{2m+4p+n-1}}}\left\| \left| x\right|
^{2p}u_{m}\right\| _{rF}.  \label{eq3g}
\end{equation}
Note that $\Delta ^{p}\left( P_{2p}f_{m}-\left| x\right|
^{2p}u_{m}\right) =0;$ thus there exists a homogeneous polynomial
$r_{m+2p}$ of degree $m+2p$ such that $\left| x\right|
^{2p}u_{m}=P_{2p}f_{m}+r_{m+2p}$ and $\Delta ^{p}r_{m+2p}=0.$ A
result proved in \cite{Rend05} (see Theorem 12 and, in particular,
equation (26), loc.\ cit.) yields the following estimate for
$f_{m}$,
\begin{equation}
\left\| f_{m}\right\| _{rF}\leq
\delta^{-1}\frac{I_{2m+n-1}}{I_{2m+2p+n-1}}\left\| \left| x\right|
^{2p}u_{m}\right\|_{rF}   \label{eqqq}
\end{equation}
where $\delta >0 $ is  a constant independent of $ m . $
So we obtain from (\ref{eq3g})
\begin{equation*}
\left\| \Delta ^{p}\left( P_{2p}f_{m}\right) \right\| _{rF}\geq \delta
 e_{p,m}
\sqrt{\frac{I_{2m+n-1}}{I_{2m+4p+n-1}}}\frac{I_{2m+2p+n-1}}{I_{2m+n-1}}
\left\| f_{m}\right\| .
\end{equation*}
It is easy to see from Proposition \ref{Iest} that there is a
constant $C'$, depending only on $p$ and $n$, such that
$$
\sqrt{\frac{I_{2m+n-1}}{I_{2m+4p+n-1}}}\frac{I_{2m+2p+n-1}}{I_{2m+n-1}}\geq
C'
$$
for all natural numbers $ m . $
This proves the estimate \eqref{eqbas1}.

The estimate \eqref{eqbas2} follows immediately from \eqref{eqbas1}
and Proposition \ref{Pest}, since there is a constant $C''>0$ such
that, for $m\geq 0$,
$$
\frac{e_{p,m}}{(m+1)^p}\geq C''.
$$
\end{proof}

We note, by Proposition \ref{Iest}, that the constant $e_{p,m}$ can
be expressed by means of the integrals $I_{m}$,
\begin{equation}\label{Iande}
\frac{I_{2m+2p+n-1}}{I_{2m+n-1}}=\frac{1}{2^{p}}\left(
2p-2+n+2m\right) ....\left( n+2m\right) =\frac{1}{2^{2p}p!}\,
e_{p,m}.
\end{equation}
We also record here the following corollary of Theorem
\ref{ThmHomA}, which will be used to prove Theorem \ref{ThmMain2}.

\begin{corollary}
\label{CorHomA}Suppose that $P_{2p}\left( x\right) $ is homogeneous
of degree $2p$ and $P_{2p}\left( x\right) \geq \delta\left| x\right|
^{2p}$ for all $x\in \mathbb{R}^{n}.$ Then there exists a constant
$D$ such that, for each homogeneous polynomial $f_{m}$ of degree
$m\in \mathbb{N}_{0},$
\begin{equation}
\left\| \Delta ^{p}\left( P_{2p}\cdot f_{m}\right) \right\|
_{rF}\geq Dm^{p}\left\| f_{m}\right\| _{rF}  \label{eqbas3}
\end{equation}
\end{corollary}

\subsection{A comment on the difference between the real and complex Fischer
 norms}
\label{comments} The following estimates for the complex Fischer
norm were proved in \cite{KhSh92},
\begin{equation}\label{CFest}
\|\Delta_\bC^p(\Sigma^p \cdot q_m)\|_F\geq C\sqrt{m^p}\|\Sigma^p\cdot
q_m\|,\quad \|\Delta_\bC^p(\Sigma^p \cdot q_m)\|_F\geq
Cm^{p}\|q_m\|_F,
\end{equation}
where $\Sigma$ is given by \eqref{standard}. These estimates lead to
the fact, mentioned in Section \ref{mainresults}, that the mixed
Cauchy problem, for $L$ with principal part $\Delta_\bC^p$, with
divisor $\Sigma$ is well posed for $k_0=3p/2$ (\cite{EbSh96}). We
note that the second estimate in \eqref{CFest} is analogous to the
estimate for the real norm in Proposition \ref{standest}. The first,
however, does not hold for the real norm in view of the following
result.
\begin{proposition}
\label{dimone}
Assume  that  $n>1 . $
Suppose that for an integer $l\geq 0$ the following estimate holds
for homogeneous polynomials $f_m$ of degree $m$:
\begin{equation}
\left\| \Delta ^{p}\left( \left| x\right| ^{2p}\cdot f_{m}\right)
\right\| _{rF}\geq C\sqrt{m^{l}}\left\| \left| x\right| ^{2p}\cdot
f_{m}\right\| _{rF} \label{eqschoen}
\end{equation}
Then $l=0.$
\end{proposition}

\begin{proof}
Suppose that \eqref{eqschoen} holds. Since $n>1 $ we may take for $f_{m}$ a
homogeneous harmonic polynomial $Y_{m}\neq 0 $ of degree  $m . $  Recall that
 $Y_{m}$ is an
eigenvector of $F_{2p}=\Delta ^{p}\left( \left| x\right| ^{2p}\cdot
\right) $ and
\begin{equation*}
\Delta ^{p}\left( \left| x\right| ^{2p}\cdot Y_{m}\right)
=d_{p}(p,m)Y_{m}
\end{equation*}
where $d_{p}(p,m)$ is given by \eqref{eigen2}. Further
\begin{equation*}
\left\| \left| x\right| ^{2p}\cdot Y_{m}\right\| _{rF}=\sqrt{\frac{
I_{4p+2m+n-1}}{I_{2m+n-1}}}\left\| Y_{m}\right\| _{rF}
\end{equation*}
Hence (\ref{eqschoen}) implies that
\begin{equation*}
| d_{p}(p,m)| \geq C\sqrt{m^{l}}\sqrt{\frac{I_{4p+2m+n-1}}{
I_{2m+n-1}}}
\end{equation*}
But $| d_{p}(p,m)| \leq Am^{p}$ and
\begin{equation*}
\frac{I_{4p+2m+n-1}}{I_{2m+n-1}}\geq m^{2p}.
\end{equation*}
So we obtain that $C\sqrt{m^{l}}\leq A,$ which implies that $l=0.$
\end{proof}

Proposition
\ref{dimone} is not true for $n=1 . $ Indeed, in this case, the set
$\mathcal{P}_{m}\left( \mathbb{R}^{n}\right) $ consists of multiples
of the polynomial  $x ^m $ and one has
\begin{equation*}
{ d^{2p} \over dx^{2p} } (x^{2p}\cdot  x^m )
= (m+2p)(m+2p-1 ) ...(m+1) x^m
 .
\end{equation*}

\section{The mixed Cauchy problem for linear partial
differential operators
 in
$\bR^n$}\label{proof}

As above, we let $B_{R}:=\left\{ x\in \mathbb{R}^{n}:\left| x\right|
<R\right\} $ denote the open unit ball in $\mathbb{R}^{n}$ and
 $A\left( B_{R}\right) $ the algebra of all infinitely
differentiable functions $f:B_{R}\rightarrow \mathbb{C}$ such that
the homogeneous Taylor series $\sum_{m=0}^{\infty }f_{m}\left(
x\right) $ converges absolutely and uniformly to $f$ on compact
subsets of $B_R$, where $f_{m}$ are the homogeneous polynomials of
degree $m$ defined by the Taylor series of $f$. Introducing polar
coordinates $x=r\theta $ with $r\geq 0$ and $\theta \in
\mathbb{S}^{n}=\left\{ x\in \mathbb{R}^{n}:\left| x\right|
=1\right\} $ one can write
\begin{equation*}
f\left( r\theta \right) =\sum_{m=0}^{\infty }r^{m}f_{m}\left( \theta
 \right)
.
\end{equation*}
If $\sum_{m=0}^{\infty }r^{m}\left| f_{m}\left( \theta \right) \right| $
converges uniformly for all $\theta \in \mathbb{S}^{n-1}$ and $0\leq r\leq
\rho <R$ then there exists a majorant $M_{\rho }$ such that
\begin{equation}
\rho ^{m}\left| f_{m}\left( \theta \right) \right| \leq M_{\rho }\text{ for
all }\theta \in \mathbb{S}^{n-1}.  \label{eqmajorant}
\end{equation}
It is easy to see that this implies
\begin{equation}
\lim \sup_{m}\max_{\theta \in \mathbb{S}^{n-1}}\sqrt[m]{\left| f_{m}\left(
\theta \right) \right| }\leq R^{-1}.  \label{hadamard}
\end{equation}
Conversely, if the estimate (\ref{hadamard}) holds for a
real-analytic function $f$ in a neighborhood of $0,$ it is easy to
see that $f\in A\left( B_{R}\right)$. We shall need the following
lemma, which follows easily from Proposition 11 in \cite{Rend05}.

\begin{proposition}
\label{convergence}Suppose that $f_{m}$ are homogeneous polynomials
of degree $m$ for $m\in \mathbb{N}_{0}.$ Then $\sum_{m=0}^{\infty
}f_{m}$ converges uniformly on compact subsets of $B_{R}$ if and only
if, for every $0<\rho<R$, there is a constant $C_\rho$ such that
\begin{equation}
 \max_{\theta \in \mathbb{S}
^{n-1}}\left| f_{m}\left( \theta \right) \right| \leq C_\rho
\rho^{-m},\quad \forall\, m\in \bN_0, \label{eqmaxnorm}
\end{equation}
if and only if, for every $0<\rho<R$, there is a constant $C_\rho$
such that
\begin{equation}
\left\| f_{m}\right\|_{rF} \leq  C_\rho \rho^{-m} \sqrt{m!},\quad
\forall\, m\in \bN_0. \label{eqFischernorm}
\end{equation}
\end{proposition}

Our main result in this section is the following. Theorem
\ref{ThmMain2} follows directly from this result in view of
Corollary \ref{CorHomA}.

\begin{theorem}
\label{abstract} Let $P_{k}$ and $Q_{k}$ be homogeneous polynomials
of degree $k$ and suppose that there exist a constant $C>0$ and an
exponent $p$ with $0<p\leq k$ such that, for all homogeneous
polynomials $q_{m}$ of degree $m\geq 0$,
\begin{equation}
\left\| Q_{k}\left( D\right) \left( P_{k}q_{m}\right) \right\|
_{rF}\geq Cm^{p}\left\| q_{m}\right\| _{rF}.  \label{eqinj}
\end{equation}
Let $0\leq k_{0}<k$ be an integer, $R>0$ a positive number,
$a_{\alpha }(x)$
functions in $A(B_{R})$ for every multi-index $\alpha \in \bN_{0}^{n}$ with
 $
|\alpha |\leq k_{0}$, and
\begin{equation}
L=Q_{k}\left( D\right) +\sum_{|\alpha |\leq k_{0}}a_{\alpha
}(x)D^{\alpha }. \label{eqdefQ}
\end{equation}
If $k_{0}<p$ then the operator $q\mapsto L\left( P_{k}q\right) $ is
a
bijection from $A\left( B_{R}\right) $ onto $A\left( B_{R}\right) $. If $
k_{0}=p$ (in particular $p<k$) then there exists $r>0$ such that
$L\left( P_{k}q\right) =f$ has a unique solution $q\in A\left(
B_{r}\right) $ for every $f\in A\left( B_{R}\right) $.
\end{theorem}

\begin{proof} Let $f=\sum_{m=0}^{\infty }f_{m}$ be a function in
 $A(B_{R})$ given in terms
of its homogeneous Taylor series as above. Consider the equation
\begin{equation}
L(P_{k}q)=f.  \label{PDE}
\end{equation}
We shall look for a solution $q$ in terms of its homogeneous Taylor series
 $
\sum_{m=0}^{\infty }q_{m}$. To prove Theorem \ref{abstract}, it
suffices to
show that the homogeneous polynomials $q_{m}$ are uniquely determined by
\eqref{PDE} and that the series $\sum_{m=0}^{\infty }q_{m}$
converges uniformly on compact subsets of $B_{r}$, for some $r>0$,
and that one can take $r=R$ if $k_{0}<p$. Let us fix $m\geq 0$ and
identify the homogeneous
part of degree $m$ in \eqref{PDE}. To this end, we expand the coefficients
 $
a_{\alpha }$ in terms of their Taylor series, $a_{\alpha
}=\sum_{m=0}^{\infty }a_{\alpha ,m}$, and obtain from \eqref{PDE}
\begin{equation}
Q_{k}(D)(P_{k}q_{m})=f_{m}-\sum_{l=0}^{k_{0}}\sum_{|\alpha
|=l}\sum_{i=0}^{m+l-k}a_{\alpha ,i}D^{\alpha }(P_{k}q_{m+l-k-i}),
\label{PDEdegm}
\end{equation}
where of course the last sum only occurs for those $l$ (if any) for which $
m+l-k\geq 0$. Note that \eqref{eqinj} implies, in particular, that $
q_{m}\mapsto Q_{k}(P_{k}q_{m})$ is injective and, hence, also
surjective as an operator from the vector space of homogeneous
polynomials (including the
zero polynomial) into itself. Since $k_{0}<k$, we conclude from
\eqref{PDEdegm} that $q_{m}$ is uniquely determined by $q_{j}$, with
$0\leq j\leq m-1$, and $f_{m}$, and that $q_{0}$ is uniquely
determined by $f_{0}$. This proves the injectivity of $q\mapsto
L\left( D\right) \left( P_{k}q\right) $.

To prove the existence of a solution to \eqref{PDE}, we must estimate the
 $\Vert \cdot \Vert _{rF}$
-norms of $q_{m}$. For the remainder of this proof, we shall only
deal with the norm $\Vert \cdot \Vert _{rF}$ and, for simplicity of
notation, shall denote this norm simply by $\Vert \cdot \Vert $. The
inequality \eqref{eqinj} implies that, for $m\geq 1$,
\begin{equation}
\begin{aligned} \|q_m\|\leq &\, C^{-1}m^{-p}\|Q_k(D)(P_kq_m)\|\\ \leq&\,
C^{-1}m^{-p}\left(\|f_m\|+\sum_{l=0}^{k_0}\sum_{|\alpha|=l}
\sum_{i=0}^{m+l-k}\|a_{\alpha,i}D^\alpha(P_kq_{m+l-k-i})\|\right).
\end{aligned}  \label{firstest}
\end{equation}
Let us fix a radius $r>0$ with $r\leq R$. To show that
$q=\sum_{m=0}^{\infty }q_{m}$ converges to a function in $A(B_{r})$,
it suffices, in view of Proposition \ref{convergence}, to show that
for every $0<\rho <r$ there is a constant $B=B_{\rho }$ such that
\begin{equation}
\Vert q_{j}\Vert \leq B\rho ^{-j}\sqrt{j!}  \label{goal}
\end{equation}
for all $j=1,2\ldots $. Let us choose two radii $\rho $ and $\sigma $ with
 $
0<\rho <\sigma <r.$ Since $f\in A(B_{R})$, there is, in view of
Proposition \ref{convergence}, a constant $D=D_{\rho }$ such that
\begin{equation}
\Vert f_{m}\Vert \leq D\rho ^{-m}\sqrt{m!},\quad m=0,1,2\ldots .
\label{fmest}
\end{equation}
Moreover, since $a_{\alpha }\in A(B_{R})$, there is, in view of
Proposition \ref{convergence}, a constant $E_{\alpha }=E_{\alpha
,\sigma }$ such that
\begin{equation}
\max_{\theta \in \bS^{n-1}}|a_{\alpha ,i}(\theta )|\leq E_{\alpha
}\sigma ^{-i},\quad i=0,1,2\ldots .  \label{aalphaest}
\end{equation}
Proposition \ref{Pest} (with $P_{k}$ replaced by $a_{\alpha ,i}$ and using that
 $  |\alpha| =l $),
shows that
\begin{equation}
\Vert a_{\alpha ,i}D^{\alpha }(P_{k}q_{m+l-k-i})\Vert \leq \sqrt{\frac{
I_{2m+n-1}}{I_{2m-2i+n-1}}}E_{\alpha }\sigma ^{-m}\Vert D^{\alpha
}(P_{k}q_{m+l-k-i})\Vert .
\label{NH1}
\end{equation}
Proposition \ref{PropDeriv} applied to $\Vert D^{\alpha
}(P_{k}q_{m+l-k-i})\Vert $ and Proposition \ref{Pest} applied to
$\Vert P_{k}q_{m+l-k-i}\Vert $, denoting the constant
$D_{k,n}M_{\bR}$ in the latter simply by $M$, yield
\begin{equation}
\Vert D^{\alpha }(P_{k}q_{m+l-k-i})\Vert \leq M2^{l/2}(m+l-i)^{l/2}\sqrt{
\frac{I_{2m+2l-2i+n-1}}{I_{2m+2l-2k-2i+n-1}}}\Vert q_{m+l-k-i}\Vert
.
\label{NH2}
\end{equation}
Let us define $n^{\ast }$ to be the natural number $\frac{1}{2}
n-1$ for even $n$ and $n^{\ast }=\left( n-1\right) /2$ for odd $n.$
It is easy to see that
\begin{equation*}
\frac{I_{2m+2l-2i+n-1}}{I_{2m+2l-2k-2i+n-1}}\leq \frac{\left(
m+l-i+n^{\ast }\right) !}{\left( m+l-k-i+n^{\ast }\right) !} .
\end{equation*}

If we set $A:=2+n^{\ast }$, then $t+n^{\ast }\leq A\cdot t$ for all
natural numbers $t\geq 1$. From this we obtain the estimate
\begin{equation*}
\frac{\left( m+l-i+n^{\ast }\right) !}{\left( m+l-k-i+n^{\ast }\right) !}%
\leq A^{k}\frac{\left( m+l-i\right) !}{\left( m+l-k-i\right) !}.
\end{equation*}
Since $m+l-i\leq m+l-i+s$ for $s=1,...,l$ we finally obtain
\begin{equation}
\frac{I_{2m+2l-2i+n-1}}{I_{2m+2l-2k-2i+n-1}}(m+l-i)^{l}\leq A^{k}\frac{%
\left( m-i+2l\right) !}{\left( m+l-k-i\right) !}.
\label{NH3}
\end{equation}
Now we conclude from \eqref{firstest} in
combination with \eqref{fmest} and \eqref{NH1}, \eqref{NH2}, \eqref{NH3}  that
\begin{equation}
\Vert q_{m}\Vert \leq C^{-1}m^{-p}D\rho ^{-m}\sqrt{m!}+S_{m}
\label{basicest3}
\end{equation}
where we define
\begin{eqnarray*}
S_{m} &:&=C^{-1}m^{-p}\sum_{l=0}^{k_{0}}\sum_{|\alpha
|=l}\sum_{i=0}^{m+l-k}E_{\alpha }\sigma ^{-i}A^{k/2}M2^{l/2}\cdot \\
&&\sqrt{\frac{\left( m+n^{\ast }\right) !}{\left( m-i+n^{\ast }\right) !}
\frac{\left( m-i+2l\right) !}{\left( m+l-k-i\right) !}}\Vert
q_{m+l-k-i}\Vert
\end{eqnarray*}
If we denote by $E_{l}$ the number
\begin{equation*}
E_{l}:=A^{k/2}M2^{l/2}\binom{n+l-1}{l}\cdot \max_{|\alpha
|=l}E_{\alpha }
\end{equation*}
and observe that the number of multi-indices $\alpha \in \bN_{0}^{n}$
for which $|\alpha |=l$ is
\begin{equation*}
\binom{n+l-1}{l},
\end{equation*}
then we get
\begin{equation*}
S_{m}\leq \frac{C^{-1}}{m^{p}}\sum_{l=0}^{k_{0}}\sum_{i=0}^{m+l-k}E_{l}
\sigma ^{-i}\sqrt{\frac{\left( m+n^{\ast }\right) !\left( m-i+2l\right) !}{
\left( m-i+n^{\ast }\right) !\left( m+l-k-i\right) !}}\Vert
q_{m+l-k-i}\Vert
.
\end{equation*}
We shall show that there exists an integer $m_{0}$ with the
following property: If there exists a constant $B$ such that
\eqref{goal} holds for all $j=0,1,\ldots ,m-1$, for any $m$ with
$m\geq m_{0}$, then \eqref{goal} holds also for $j=m.$ The existence
of such an integer $m_{0}$ clearly implies, by induction, that
\eqref{goal} holds for all $j=0,1,2,\ldots $,
which in turn completes the proof of the theorem.

So suppose that, for some
 $
B$, \eqref{goal} holds for all $j=0,1,\ldots ,m-1$. Without loss of
generality, of course, we may assume $B\geq 1$. The estimate
\eqref{basicest3} can be written as
\begin{equation}
\Vert q_{m}\Vert \leq B\rho
^{-m}\sqrt{m!}(C^{-1}DB^{-1}m^{-p}+T_{m}), \label{qmT}
\end{equation}
where we have defined
\begin{equation*}
T_{m}:=\frac{C^{-1}}{m^{p}}\sum_{l=0}^{k_{0}}\sum_{i=0}^{m+l-k}E_{l}\sigma
^{-i}\sqrt{\frac{\left( m+n^{\ast }\right) !\left( m-i+2l\right)
!}{m!\left( m-i+n^{\ast }\right) !}}\rho ^{-(l-k-i)}.
\end{equation*}
Recall that $B\geq 1$, so that $C^{-1}DB^{-1}m^{-p}\leq
C^{-1}Dm^{-p}$. We shall show that there is $m_{0}$ such that
$C^{-1}Dm^{-p}+T_{m}\leq 1$ for all $m\geq m_{0}$.  For this, it
suffices to show that there exists $0<\theta <1$ and $m_{0}$ such that
$$
T_{m}\leq \theta \text{ for all } m\geq m_{0}.$$
 We set
\begin{equation*}
N\left( \rho \right) :=C^{-1}(k_{0}+1)\max_{0\leq l\leq
k_{0}}E_{l}\cdot \max_{0\leq l\leq k_{0}}\rho ^{k-l}
\end{equation*}
and obtain
\begin{equation*}
T_{m}\leq \frac{N\left( \rho \right)
}{m^{p}}\sum_{i=0}^{m+k_{0}-k}\left( \frac{\rho }{\sigma }\right)
^{i}\sqrt{\frac{\left( m+n^{\ast }\right) !\left( m-i+2k_{0}\right)
!}{m!\left( m-i+n^{\ast }\right) !}}.
\end{equation*}
First suppose that $2k_{0}\geq n^{\ast }$. Then
\begin{equation*}
T_{m}\leq \frac{N\left( \rho \right) }{m^{p}}\left( m+n^{\ast
}\right) ^{n^{\ast }/2}\left( m+2k_{0}\right) ^{\left(
2k_{0}-n^{\ast }\right) /2}\sum_{i=0}^{\infty }\left( \frac{\rho
}{\sigma }\right) ^{i}.
\end{equation*}
If $k_{0}<p$, then the right hand side clearly converges to zero as
$m\rightarrow \infty$. If $ k_{0}=p$, then we can make $T_m\leq
\theta<1$ for all sufficiently large $m$ by making $\rho<r$
sufficiently small. (Observe that $N(\rho)\to 0$ as $\rho\to 0$,
since $k_0<k$.)

Now suppose that $2k_{0}<n^{\ast }.$ Clearly we obtain
\begin{equation}
T_{m}\leq \frac{N\left( \rho \right) }{m^{p}}\left( m+n^{\ast
}\right) ^{n^{\ast }/2}\sum_{i=0}^{m+k_{0}-k}\left( \frac{\rho
}{\sigma }\right) ^{i} \frac{1}{\sqrt{\left( m-i+2k_{0}+1\right)
\ldots \left( m-i+n^{\ast }\right) }}. \label{eqnewsum}
\end{equation}
In order to estimate the latter sum, we fix a number $\delta$ with
$0<\delta<1$. We consider only $m$ with $m>\delta ^{-1}.$ The
estimate $m-i\geq \left( 1-\delta \right) m$ holds for all
 $i\leq \delta m.$ We split the sum in (\ref{eqnewsum}) into two sums $
I_{1}+I_{2}$, the first one containing only indices with $i\leq \left[
 \delta m
\right] $ and the second one indices $i$ with $i>\left[ \delta
m\right] .$ We estimate
\begin{equation*}
I_{1}\leq \sum_{i=0}^{\left[ \delta m\right] }\left( \frac{\rho }{\sigma
 }%
\right) ^{i}\frac{1}{\left(\left( 1-\delta \right) m\right)^{\left(
n^{\ast }-2k_{0}\right) /2}}\leq K\frac{1}{ m^{\left(
n^{\ast }-2k_{0}\right) /2}}\frac{1}{
1-\frac{\rho }{\sigma }}
\end{equation*}
where
\begin{equation*}
K:=\frac{1}{\left( 1-\delta \right) ^{\left( n^{\ast }-2k_{0}\right)
/2}} .
\end{equation*}
For the second sum $I_{2}$, we use the estimate
\begin{equation*}
I_{2}\leq \sum_{i=\left[ \delta m\right] }^{m+k_{0}-k}
\left( \frac{\rho
 }{
\sigma }\right) ^{i}\leq \left( \frac{\rho }{\sigma }\right)
^{\left[ \delta m\right] }\frac{1}{1-\frac{\rho }{\sigma }}
\end{equation*}
and, hence, obtain
\begin{equation*}
T_{m}\leq \frac{N\left( \rho \right) }{1-\frac{\rho }{\sigma
}}\frac{\left( m+n^{\ast }\right) ^{n^{\ast }/2}}{m^{p}}\left(
\frac{K}{m^{\left( n^{\ast }-2k_{0}\right) /2}}+\left( \frac{\rho
}{\sigma }\right) ^{[\delta m]}\right).
\end{equation*}
As before it is easy to see that $T_{m}$ converges to $0$ for
$k_{0}<p$, since $m^{s}\left( \frac{\rho }{\sigma }\right) ^{[\delta
m]}$ converges to $0$ for any integer $s$ (recall that $\rho /\sigma
<1).$ If $k_{0}=p$, we use, as above, the fact that $N\left( \rho
\right) \rightarrow $ as $\rho \rightarrow 0$ to conclude that, if
$r$ is sufficiently small, then $T_m\leq \theta<1$ for large $m$.
This completes the proof of Theorem \ref{abstract}.
\end{proof}

\end{document}